\documentclass[12pt]{amsart}
\usepackage{hyperref}
\hypersetup{nesting=true,debug=true,naturalnames=true}
\usepackage{graphicx,amssymb,upref}

\let\<\langle
\let\>\rangle

\let\uml\"

\usepackage{amsmath,amsthm,amscd,amsfonts,amssymb,enumerate,pgf,tikz,mathrsfs,caption}
\theoremstyle{plain}
\newtheorem{theorem}{Theorem}[section]

\newtheorem{prop}[theorem]{Proposition}

\theoremstyle{definition}
\newtheorem{definition}[theorem]{Definition}
\newtheorem{example}[theorem]{Example}

\newtheorem{corollary}[theorem]{Corollary}

\theoremstyle{remark}
\newtheorem{remark}[theorem]{Remark}

\numberwithin{equation}{section}

\newcommand{\bN}{\mathbb{N}}

\newcommand{\Dom}{\mathrm{Dom}}

\newcommand{\cA}{\mathcal{A}}

\newcommand{\af}{\ast_f}
\begin{document}

\title{SOME MONOCHROMATIC PATTERNS IN NATURAL NUMBERS}

%    Information for first author
\author{Arpita Ghosh and Surojit Ghosh}

\curraddr{Department of Mathematics, University of Haifa, Israel}
\email{arpi.arpi16@gmail.com}

\curraddr{Department of Mathematics, IIT Roorkee, Roorkee}

\email{surojit.ghosh@ma.iitr.ac.in; surojitghosh89@gmail.com}

\subjclass[2020]{05D10, 05C55, 22A15}
%\date{January 1, 2001 and, in revised form, June 22, 2001.}

\keywords{Ramsey Theory, Monochromatic patterns, The set of sums of two squares, Algebra in the Stone-\v{C}ech compactification}

\maketitle
\begin{abstract}
The set of sums of two squares plays a significant role in elementary number theory. In this article, we establish the existence of several rich monochromatic configurations in the natural numbers by exploiting algebraic structures induced by the set of sums of two squares. The proofs rely largely on the algebraic properties arising from the induced structures on the Stone-\v{C}ech compactification of the natural numbers.

%In \cite{Nasso22}, Di Nasso investigated various infinite monochromatic %patterns in integers using operations induced by affine maps and posed the %question of whether different sets of infinite monochromatic configurations %could be found in natural numbers under the structure induced by the set of %sums of two squares. This article addresses Nasso's question and, in doing %so, establishes several classical Ramsey-type theorems in this novel setting.
\end{abstract}

\section{Introduction}
The search for monochromatic patterns within colorings of algebraic structures has been a longstanding theme in Ramsey theory. Among the fundamental results in this area, the celebrated \emph{van der Waerden theorem} guarantees the existence of monochromatic arithmetic progressions for any finite partition of the natural numbers. This theorem illustrates that the family of arithmetic progressions is \emph{partition regular}. Notably, the van der Waerden theorem is a finitary Ramsey-type result relying solely on the additive structure of $\mathbb{N}$. Its infinitary counterpart, \emph{Hindman's theorem}, asserts that for any finite coloring of $\mathbb{N}$, there exists an injective sequence whose finite sums are all monochromatic.

A natural direction in Ramsey theory is the transition from additive to multiplicative structures. One such transition can be achieved via the transformation $n \mapsto p^n$ for a prime $p$, leading to the \emph{multiplicative van der Waerden theorem}, which states that for any finite partition of $\mathbb{N}$, a monochromatic geometric progression of arbitrary length exists.

Patterns that intertwine the additive and multiplicative structures of natural numbers present greater challenges. The pioneering work of Bergelson \cite{Ber10} and Hindman \cite{Hind11} independently established the monochromaticity of the configuration ${a, b, c, d}$ where $a + b = c \cdot d$. Further progress was made by Moreira \cite{Joel17}, who demonstrated that the pattern ${a, a + b, a \cdot b}$ is monochromatic. Subsequently, Moreira, Barrett, and Lupini \cite{BLM21} extended this result by proving the monochromaticity of ${a, a + b, a + b + a \cdot b}$. In \cite{Saha18}, Sahasrabudhe proved rich exponential patterns in $\mathbb{N}$, showing in particular that the sets $\{a,b, a^b \}$ and $\{a,b, ab, a^b\}$ are monochromatic where $b >1$. Subsequently, Goswami in \cite{G25} established the result that $\{a, b, ab, (a+1)b\}$ is monochromatic. However, the question of whether the set $\{a, b, a + b, a \cdot b\}$ is always monochromatic remains an open problem.

Recently, in \cite{Nasso22}, Nasso found many new additive and multiplicative monochromatic patterns by introducing the notions of symmetric polynomials. Polynomial extension of some important classical results, particularly the polynomial van der Waerden Theorem and the polynomial Deuber Theorem for symmetric polynomials, have been explored by Chakraborty and Goswami in \cite{CG24}.

Given these advancements, it is natural to examine the Ramsey-theoretic properties of well-structured subsets of natural numbers. Define $S(m,n)$ as the set:
$$
S(m,n):=\{x_1^m+\cdots + x_n^m: x_1, \cdots, x_m \in \mathbb{N}_0\}.
$$
From elementary number theory, it is well known that $S(2,4) = \mathbb{N}_0$, while $S(2,3)$ consists of natural numbers that are not of the form $4^a(8b+7)$ for some $a, b \in \mathbb{N}_0$. The set $S(2,3)$ lacks multiplicative closure, whereas $S(2,4)$ is multiplicatively closed but trivial in the sense that it encompasses all natural numbers. Thus, our focus is on $\Sigma = S(2,2)$, the set of all natural numbers whose prime factors of the form $4k+3$ (if any) have even exponents. 

In this work, we introduce a new algebraic structure on natural numbers derived from the set of sums of two squares, $\Sigma$, and, employing the Stone-\v{C}ech compactification, establish several classical results, including Hindman's theorem, Deuber's theorem, Brauer's theorem, Milliken--Taylor's theorem, the existence of geo-arithmetic progressions, and a polynomial variant of van der Waerden's theorem.

\section{A new structure and  Stone-\v{C}ech compactification} 

Let \(\Sigma = \{a^2 + b^2 : a, b \in \mathbb{N}_0\} \subseteq \mathbb{N}_0\) denote the set of non-negative integers expressible as sums of two squares. Under standard multiplication, \(\Sigma\) forms a \textbf{commutative semigroup}, as the product of two elements in \(\Sigma\) remains in \(\Sigma\) by the identity:
\[
(a^2 + b^2)(c^2 + d^2) = (ad - bc)^2 + (ac + bd)^2.
\]
We order \(\Sigma\) increasingly as:
\[
\Sigma = \{s_0 < s_1 < s_2 < s_3 < \cdots\},
\]
where \(s_0 = 0\), \(s_1 = 1\), \(s_2 = 2\), \(s_3 = 4\), \(s_4 = 5\), \(s_5 = 8\), \(s_6 = 9\), \(s_7 = 10\), \(s_8 = 13\), \(s_9 = 16\), \(s_{10} = 17\), \(s_{11} = 18\), \(s_{12} = 20\), \(s_{13} = 25\), \(s_{14} = 26\), \(s_{15} = 29\), \(s_{16} = 32\), and so forth.

For \(s \in \Sigma\), define the \textbf{predecessor set} \(\Sigma_{<s} = \{y \in \Sigma : y < s\}\). Let \(g: \Sigma \to \mathbb{N}_0\) be the function:
\[
g(s) = \left|\Sigma_{<s}\right|,
\]
which counts the number of elements in \(\Sigma\) strictly less than \(s\). This \(g\) is a bijection, with inverse \(f: \mathbb{N}_0 \to \Sigma\) given by \(f(n) = s_n\), the \(n\)-th term in the ordered sequence of \(\Sigma\).

\begin{definition}[Induced Operation]\label{oper}
The operation \(\ast_f\) on \(\mathbb{N}_0\) is defined as:
\[
m \ast_f n = g(f(m) \cdot f(n)) = \left|\Sigma_{<s_m \cdot s_n}\right|.
\]
\end{definition}

\begin{example}
Let \(m = 2\) and \(n = 5\). Then \(s_m = 2\) and \(s_n = 8\). Their product is \(s_m \cdot s_n = 16 = s_9\). The set \(\Sigma_{<16}\) contains the first 9 elements of \(\Sigma\) (i.e., \(s_0, s_1, \ldots, s_8\)), so \(\left|\Sigma_{<16}\right| = 9\). Hence, \(2 \ast_f 5 = 9\).
\end{example}

\begin{remark}

As we note that $\Sigma$ is multiplicatively closed, the product \( s_i \cdot s_j \) is an element of \( \Sigma \), specifically:
\[
s_i \cdot s_j = s_k
\]
for some \( k \in \mathbb{N} \cup \{0\} \). The exact value of \( k \) depends on the ordering (based on the numerical value) of \( \Sigma \), which is not uniquely determined by \( i \) and \( j \) alone.
\end{remark}
\begin{prop}
\((\mathbb{N}_0, \ast_f)\) is a commutative semigroup with identity element \(1\), and \(f: (\mathbb{N}_0, \ast_f) \to (\Sigma, \cdot)\) is a semigroup isomorphism.
\end{prop}

\begin{proof}
The operation \(\ast_f\) inherits associativity and commutativity from the semigroup \((\Sigma, \cdot)\) via the bijection \(f\). Explicitly:
\begin{itemize}
\item  For all \(n \in \mathbb{N}_0\),
\[
1 \ast_f n = g(f(1) \cdot f(n)) = g(s_1 \cdot s_n) = g(s_n) = \left|\Sigma_{<s_n}\right| = n,
\]
since \(s_1 = 1\) is the multiplicative identity in \(\Sigma\).

\item  By definition, 
\[
f(m \ast_f n) = f(g(s_m \cdot s_n)) = s_{g(s_m \cdot s_n)} = s_m \cdot s_n = f(m) \cdot f(n),
\]
confirming \(f\) is a homomorphism. As \(f\) is bijective, it is an isomorphism.
\end{itemize}
\end{proof}

\noindent Define the \emph{power} \(x^{(n)}\) in \((\mathbb{N}_0, \ast_f)\) recursively by:
\begin{itemize}
\item \(x^{(0)} = 1\),
\item \(x^{(n)} = x \ast_f x \ast_f \cdots \ast_f x = \left|\Sigma_{<s_x^n}\right|\).
\end{itemize}

We briefly recall the algebraic structure of the \emph{Stone-Čech compactification} \(\beta S\) for a discrete semigroup \((S, \cdot)\). The elements of \(\beta S\) are the \emph{ultrafilters} on \(S\). By identifying the principal ultrafilters with the points of \(S\), we may regard \(S\) as a subset of \(\beta S\), i.e., \(S \subseteq \beta S\). 

For a subset \(A \subseteq S\), define:
\[
\overline{A} = \{p \in \beta S \mid A \in p\}.
\]
The collection \(\{\overline{A} \mid A \subseteq S\}\) forms a basis for a topology on \(\beta S\). The semigroup operation \(\cdot\) on \(S\) can be extended to \(\beta S\), making \((\beta S, \cdot)\) a compact right topological semigroup. This means:
\begin{itemize}
    \item For any \(p \in \beta S\), the right translation map \(\rho_p: \beta S \to \beta S\) defined by \(\rho_p(q) = q \cdot p\) is continuous.
    \item The semigroup \(S\) is contained in the topological center of \(\beta S\), meaning that for any \(x \in S\), the left translation map \(\lambda_x: \beta S \to \beta S\) defined by \(\lambda_x(q) = x \cdot q\) is continuous.
\end{itemize}

The extended operation on \(\beta S\) is characterized as follows: for \(p, q \in \beta S\) and \(A \subseteq S\),
\[
A \in p \cdot q \mbox{ if and only if }  \{x \in S \mid x^{-1} A \in q\} \in p,
\]
where \(x^{-1} A = \{y \in S \mid x \cdot y \in A\}\). A fundamental result due to Ellis states that every compact right topological semigroup contains an idempotent element. Thus, there exists an idempotent ultrafilter \(p \in \beta S\) satisfying \(p \cdot p = p\).

For a sequence \(\langle x_n \rangle_{n=1}^\infty\) in \((S, \cdot)\), the set of finite products is defined as:
\[
FP(\langle x_n \rangle_{n=1}^\infty) = \{x_{n_1} \cdot x_{n_2} \cdots x_{n_k} \mid n_1 < n_2 < \cdots < n_k\}.
\]

Idempotent ultrafilters play a crucial role in Ramsey Theory. This is exemplified by Galvin's theorem, which states that for any semigroup \((S, \cdot)\) and an idempotent ultrafilter \(p = p \cdot p\) in \((\beta S, \cdot)\), every set \(A \in p\) contains the set of finite products \(FP(\langle x_n \rangle_{n=1}^\infty)\) for some sequence \(\langle x_n \rangle_{n=1}^\infty\) in \(S\). For further details on the Stone-Čech compactification and its applications in Ramsey Theory, see \cite{HS12} and \cite{Tod10}.

Since \((\mathbb{N}_0, \af)\) is a discrete commutative semigroup, its Stone-Čech compactification \((\mathbb{N}_0, \af)\) gives rise to numerous Ramsey-type results in \(\mathbb{N}_0\) induced by the operation \(\af\).

\section{Study on Monochromatic Configurations}
In this section, we introduce new monochromatic configurations within the semigroup $(\mathbb{N}, \ast_f)$, which encodes information about sums of two squares. We begin with the following definition:

For any infinite sequence of natural numbers $\langle x_n \rangle_{n=1}^{\infty}$, the corresponding set of finite sums is denoted and defined as
$$FS(\langle x_n\rangle_{n=1}^{\infty}) = \{ x_{n_1} + x_{n_2} + \cdots + {x_{n_k}} : n_1 < n_2 < \cdots < n_k \}.$$

A fundamental result in arithmetic Ramsey Theory guarantees the existence of infinite monochromatic patterns within such sets of finite sums. This result, due to Hindman, is stated as follows:

\begin{theorem}[Hindman's Finite Sum Theorem, \cite{Hind74}]
Let $r \geq 1$. For every $r$-coloring of the natural numbers, $\mathbb{N} = C_1 \cup C_2 \cup \cdots \cup C_r$, there exists a color $C_i$ and a sequence $\langle x_n \rangle_{n=1}^{\infty}$ in $\mathbb{N}$ such that $FS(\langle x_n \rangle_{n=1}^{\infty}) \subseteq C_i.$

\end{theorem}

A similar result holds when considering finite products instead of finite sums, and more generally, in any semigroup \cite[Theorem 5.8]{HS12}. Analogously, for the semigroup $(\mathbb{N}_0, \ast_f)$, we define the set of finite $\ast_f$-operations as follows:
\[\mathfrak{FP}_f(\langle x_n \rangle_{n=1}^{\infty})= \{x_{n_1} \ast_f x_{n_2} \ast_f \cdots \ast_f x_{n_k} : n_1 < n_2 < \cdots < n_k\}
\]

Equivalently,
\[\mathfrak{FP}_f(\langle x_n \rangle_{n=1}^{\infty})= \{|\Sigma_{< s_{x_{n_1}}  s_{x_{n_2}}  \cdots s_{x_{n_k}}}| : n_1 < n_2 < \cdots < n_k\}
\]

where $\Sigma$ denotes the set of sums of two squares.

\begin{theorem}
For any $r \geq 1$ and any finite $r$-coloring of the nonnegative integers, $\mathbb{N}_0 = C_1 \cup C_2 \cup \cdots \cup C_r$, there exist a color $C_i$ and a sequence $\langle x_n \rangle_{n=1}^{\infty}$ in $\mathbb{N}$ such that

$$\mathfrak{FP}_f(\langle x_n \rangle_{n=1}^{\infty}) \subseteq C_i$$

In other words,
\[
\{|\Sigma_{< s_{x_{n_1}}  s_{x_{n_2}}  \cdots s_{x_{n_k}}}| : n_1 < n_2 < \cdots < n_k\} \subseteq C_i,
\]

\end{theorem}

\begin{proof}
The proof is actually verbatim to the proof of the Hindman's Theorem \cite[Theorem 5.8]{HS12}, replacing the usual addition with the associative operation $\ast_f$.
\end{proof}

A cornerstone result in arithmetic Ramsey Theory is Van der Waerden's Theorem (1927), which guarantees the existence of arbitrarily long monochromatic arithmetic progressions. A year later, Brauer strengthened this result by showing that one can also ensure that the common difference belongs to the same color class as the elements of the progression.

\begin{theorem}[van der Waerden's Theorem, \cite{Van27}]
For every finite coloring $\mathbb{N} = C_1 \cup C_2 \cup \cdots \cup C_r$ and for every $L \in \mathbb{N}$, there exists a monochromatic arithmetic progression of length $L$; that is, there exist a color $C_i$ and elements $a, b \in \mathbb{N}$ such that
$a, a+b, a+2b, \cdots , a+Lb \in C_i.$
\end{theorem}

\begin{theorem}[Brauer's Theorem, \cite{Brauer28}]
For every finite coloring $\mathbb{N} = C_1 \cup C_2 \cup \cdots \cup C_r$ and for every $k \in \mathbb{N}$, there exist a color $C_i$ and elements $a, b \in \mathbb{N}$ such that
 $\{a, b, a+b, a+2b, \cdots , a+kb\} \subseteq C_i.$
\end{theorem}

We establish an analogous version of Brauer's Theorem in our setting:

\begin{theorem} \label{Brauer}
For every finite coloring $\mathbb{N}_0 = C_1 \cup C_2 \cup \cdots \cup C_r$ and for every $k \in \mathbb{N}$, there exist a color $C_i$ and elements $x, z \in \mathbb{N}$ such that
$$\{x, z, |\Sigma_{<{s_x^j}s_z}|: j=1,2, \cdots, k\} \subseteq C_i.$$
\end{theorem}

After a few decades, Deuber established a result concerning the generalized partition regularity of homogeneous systems of linear Diophantine equations. In particular, he demonstrated the partition regularity of the so-called $(m,p,c)$-sets.

\begin{theorem}[Deuber's Theorem, \cite{Deuber74}]
\label{deuber}
For every $m, p, c \in \mathbb{N}$ and for every finite coloring $\mathbb{N}= C_1 \cup C_2  \cup \cdots \cup C_r$, there exists a monochromatic $(m,p,c)$-set. That is, there exists a color $C_i$ and elements $a_0, a_1, \dots, a_m \in C_i$ such that
 \[ca_j + \sum_{s=0}^{j-1} n_s a_s \in C_i \mbox{ for every } j \in \{1,2, \cdots, m \} \mbox{ and for all } n_0, n_1, \cdots , n_{j-1} \in \{-p, \cdots, p\}.\] 
\end{theorem}

We present an analogous version of Deuber's theorem in our context as follows:

\begin{theorem}
\label{newDeuber}
Let $m, p, r \in \mathbb{N}$. For every $r$-coloring $\mathbb{N}_0 = C_1 \cup C_2  \cup \cdots \cup C_r$, there exists a color $C_i$ and elements $x_0, x_1, \dots, x_m \in C_i$ such that
 \[
 \{\left|\Sigma_{<{s_{x_0}^{n_0}}{s_{x_1}^{n_1}} \cdots {s_{x_{j-1}}^{n_{j-1}}} {s_{x_j}}}\right| : n_0, n_1, \cdots, n_{j-1} \in \{0,1, \cdots, p\} \mbox{ and } j=1,2, \cdots, m \}  \subseteq  C_i
 .\]
\end{theorem}
\begin{proof}
We prove this theorem using a generalization of Deuber's theorem for commutative semirings, recently established by V. Bergelson, J.H. Johnson, and J. Moreira in \cite[Corollary 3.7]{BJM17}. It states the following:

Let $(S, \ast)$ be a commutative semigroup, and for each $j=1,2, \dots, m$, let $\mathfrak{F}_j$ be a finite set of endomorphisms $\psi : S^j \to S$. Then, for every $r$-coloring $S=C_1 \cup C_2  \cup \cdots \cup C_r$, there exist a color $C_i$ and elements $x_0, x_1, \dots, x_m$ distinct from the identity, such that $x_0 \in C_i$ and
 \[
 \psi(x_0, x_1, \cdots, x_{j-1}) \ast x_j \in C_i \mbox{ for every } j=1,2, \cdots, m \mbox{ and for every } \psi \in \mathfrak{F}_j.\]
This statement generalizes Theorem \ref{deuber} for the case $c=1$. We apply the Bergelson-Johnson-Moreira result with $(S, \ast) = (\mathbb{N}_0, \ast_f)$. For each $j$-tuple $\bar{n}=(n_0, n_1, \dots, n_{j-1}) \in (\mathbb{N}_0)^j$, define 
\[
\psi_{\bar{n}} : (x_0, x_1, \cdots, x_{j-1})\mapsto {x_0}^{(n_0)} \ast_f  {x_1}^{(n_1)} \ast_f \cdots \ast_f  {x_j}^{(n_j)}.
\]
Since $(\mathbb{N}_0, \ast_f)$ is a commutative semigroup, each $\psi_{\bar{n}}$ is a semigroup homomorphism. Define 
\[
\mathfrak{F}_j= \{\psi_{\bar{n}} : \mathbb{N}^j \rightarrow \mathbb{N} : \bar{n}=(n_0, n_1, \cdots, n_{j-1}) \in \{0,1,2, \cdots, p\}^j\}
\]
which consists of homomorphisms for each $j=1,2, \dots, m$. Then, by the Bergelson-Johnson-Moreira theorem, for every finite coloring $\mathbb{N}= C_1 \cup C_2  \cup \cdots \cup C_r$, there exist a color $C_i$ and elements $x_0, x_1, \dots, x_m$ distinct from the identity such that:

\begin{enumerate}
\item $x_0 \in C_i$,
\item $\psi_{\bar{n}}(x_0, x_1, \dots, x_{j-1}) \ast_f x_j \in C_i$ for every $j=1,2, \dots, m$ and all $\bar{n} \in \{0,1,2, \dots, p\}^j$.
\end{enumerate}
From the second condition, we deduce that 
\[
\psi_{\bar{n}}(x_0,x_1, \cdots, x_{j-1}) \ast_f x_j \in C_i
\]
Equivalently, 
\[
{x_0}^{(n_0)} \ast_f  {x_1}^{(n_1)} \ast_f \cdots \ast_f  {x_j}^{(n_j)} \ast_f x_j \in C_i.
\]

Therefore, for all $(n_0, n_1, \cdots, n_{j-1}) \in \{0,1, \cdots, p\}^j \mbox{ and } j=1,2, \cdots, m$ we obtain 
$$\left|\Sigma_{<{s_{x_0}^{n_0}}{s_{x_1}^{n_1}} \cdots {s_{x_{j-1}}^{n_{j-1}}} {s_{x_j}}}\right| \in C_i.$$ 
Hence, the result follows.
\end{proof}

\begin{proof}(Proof of Brauer's Theorem \ref{Brauer})The proof follows from the Bergelson-Johnson-Moreira theorem, as applied in Theorem \ref{newDeuber}. For $j=0,1,2, \dots, k+1$, let $\psi_j : (\mathbb{N}, \ast_f) \to (\mathbb{N}, \ast_f)$ be the endomorphism defined by $\psi_j(x) = x^{(j)}$. Then, for every $r$-coloring $\mathbb{N}= C_1 \cup C_2  \cup \cdots \cup C_r$, there exist a color $C_i$ and elements $x, y \neq 1$ such that
\begin{equation}\label{deuequ}
\{x, \psi_0(x) \ast_f y, \psi_1(x) \ast_f y, \psi_2(x) \ast_f y , \cdots, \psi_{k+1}(x) \ast_f y\} \subseteq C_i .
\end{equation}

Setting $z:= x \ast_f y$, we obtain 
\[
\{x, z, x\ast_f z, x^{(2)} \ast_f z, x^{(3)} \ast_f z , \cdots, x^{(k)} \ast_f z \} \subseteq C_i.
\]
In other words,
\[
\{x, z, \left|\Sigma_{<s_x^js_z}\right|: j=1,2, \cdots, k\} \subseteq C_i.
\]
Hence the result follows.
\end{proof}

\section{Milliken--Taylor Theorem}

The Milliken--Taylor Theorem simultaneously generalizes Hindman's Finite Sums Theorem and Ramsey's Theorem. It has been widely applied in combinatorial number theory, including extensions of Szemerédi's Theorem on arithmetic progressions.

To state the theorem, we first establish notation. Let $\mathcal{P}_f(\mathbb{N})$ denote the set of all finite non-empty subsets of $\mathbb{N}$. For $m \in \mathbb{N}$, define $[\mathbb{N}]^m = \{A \subseteq \mathbb{N} : |A| = m\}$. For $F, G \in \mathcal{P}_f(\mathbb{N})$, we write $F < G$ if $\max F < \min G$.

\begin{theorem}[Milliken--Taylor Theorem \cite{Milli75, Tay76}]\label{thm:Milliken--Taylor}
  Let $r, m \in \mathbb{N}$ and let $[\mathbb{N}]^m = C_1 \cup \cdots \cup C_r$ be an $r$-coloring. There exist an injective sequence $\langle x_n \rangle_{n=1}^\infty$ in $\mathbb{N}$ and a color $C_i$ such that  
  \[
    \left\{(x_{F_1}, \dots, x_{F_m}) : F_1 < \cdots < F_m\right\} \subseteq C_i,
  \]
  where for $F = \{n_1 < \cdots < n_k\} \in \mathcal{P}_f(\mathbb{N})$, we define $x_F = \sum_{i=1}^k x_{n_i}$.
\end{theorem}

\textbf{Remark:} When $m=1$, this reduces to Hindman's Theorem. If each $F_i$ is a singleton, we recover Ramsey's Theorem.

\subsection{Ultrafilter Formulation}
The theorem admits an ultrafilter-theoretic formulation using tensor products. Let $(S_i, \cdot)$ be semigroups and $p_i \in \beta S_i$ ultrafilters. The \textit{tensor product} $\bigotimes_{i=1}^k p_i \in \beta\left(\times_{i=1}^k S_i\right)$ is defined inductively:
\begin{itemize}
  \item For $k=1$, $\bigotimes_{i=1}^1 p_i = p_1$.
  \item For $k \geq 1$, $A \in \bigotimes_{i=1}^{k+1} p_i$ iff 
  \[
    \left\{(x_1, \dots, x_k) \in \times_{i=1}^k S_i : \{x_{k+1} \in S_{k+1} : (x_1, \dots, x_{k+1}) \in A\} \in p_{k+1}\right\} \in \bigotimes_{i=1}^k p_i.
  \]
\end{itemize}

\begin{theorem}[Bergelson-Hindman-Williams \cite{BHW14}]\label{thm:BHW}
  Let $S$ be a semigroup, $m \in \mathbb{N}$, and $A \subseteq \times _{i=1}^m S$. The following are equivalent:
  \begin{enumerate}
    \item There exists a sequence $\langle x_n \rangle_{n=1}^\infty$ in $S$ with $\{(x_{F_1}, \dots, x_{F_m}) : F_1 < \cdots < F_m\} \subseteq A$.
    \item There exists an idempotent ultrafilter $p \in \beta S$ such that $A \in \bigotimes_{i=1}^m p$.
  \end{enumerate}
\end{theorem}

For a function $\phi: X \to Y$, the induced map $\phi_*: \beta X \to \beta Y$ is defined by $\phi_*(p) = \{B \subseteq Y : \phi^{-1}(B) \in p\}$. Note that $A \in p$ implies  $\phi(A) \in \phi_*(p)$, but the converse fails.

\begin{corollary}\label{cor:tensor-product}
  Let $(S, \cdot)$ be a semigroup, $m \in \mathbb{N}$, $\phi: S^m \to S$, and $B \subseteq S$. The following are equivalent:
  \begin{enumerate}
    \item There exists a sequence $\langle x_n \rangle_{n=1}^\infty$ in $S$ with $\{\phi(x_{F_1}, \dots, x_{F_m}) : F_1 < \cdots < F_m\} \subseteq B$.
    \item $B \in \phi_*\left(\bigotimes_{i=1}^m p\right)$ for some idempotent $p \in \beta S$.
  \end{enumerate}
\end{corollary}

\begin{proof}
  Apply Theorem~\ref{thm:BHW} with $A = \phi^{-1}(B)$. For details, see \cite[Corollary 5.4]{Nasso22}.
\end{proof}

\begin{theorem}\label{thm:main}
  Let $r, m \in \mathbb{N}$, $\phi: (\mathbb{N}_0)^m \to \mathbb{N}_0$, and $\mathbb{N}_0 = C_1 \cup \cdots \cup C_r$ be an $r$-coloring. There exists a sequence $\langle x_n \rangle_{n=1}^\infty$ in $\mathbb{N}$ and a color $C_i$ such that 
  \[
   \left\{\phi\left(|\Sigma_{<F_1}|, \dots, |\Sigma_{<F_m}| \right) : F_1 < \cdots < F_m \right\} \subseteq C_i,\]
  where for $F_j = \{n_{j_1} < \cdots < n_{j_{k_j}}\} \in \mathcal{P}_f(\mathbb{N})$, we define $\Sigma_{<F_j} = \Sigma_{<s_{x_{n_{j_1}}} \cdots s_{x_{n_{j_{k_j}}}}}$.
\end{theorem}

\begin{proof}
  Let $(S, \cdot) = (\mathbb{N}_0, \ast_f)$, where $\ast_f$ is defined by $x \ast_f y = |\Sigma_{<s_x s_y}|$. Let $p \in \beta \mathbb{N}_0$ be an idempotent ultrafilter. By Corollary~\ref{cor:tensor-product}, with $q = \phi_*\left(\bigotimes_{i=1}^m p\right)$, there exist $C_i \in q$ and a sequence $\langle x_n \rangle_{n=1}^{\infty}$ such that  
  \[
    \{\phi(x_{F_1}, \dots, x_{F_m}) : F_1 < \cdots < F_m\} \subseteq C_i.
  \]
  For $F_j = \{n_{j_1} < \cdots < n_{j_{k_j}}\}$, compute:
  \[
    x_{F_j} = x_{n_{j_1}} \ast_f \cdots \ast_f x_{n_{j_{k_j}}} = \left|\Sigma_{<s_{x_{n_{j_1}}} \cdots s_{x_{n_{j_{k_j}}}}}\right| = |\Sigma_{<F_j}|.
  \]
  Substituting into $\phi$ yields the result.\end{proof}

\section{Application of Hales-Jewett theorem}
We start this section with definition of \textit{partial semigroup}:

\begin{definition}
A \emph{partial semigroup} is a triple $(S, X \subseteq S \times S, \ast)$ of a set $S$, a subset $X \subseteq S \times S$ and an operation  $\ast$ 
defined on $X$ satisfying
\[
(x \ast y) \ast z = x \ast (y \ast z)\hspace{1cm} \forall x,y,z \in G
\]
 in the sense that if either side 
is defined, so is the other and they are equal.

\end{definition}

 Let $\cA$ be a nonempty, finite set of alphabet and $v \notin \cA$ be a symbol, call it a variable. A located word in the alphabet $\cA$ is finitely supported function $b : \Dom(b) \rightarrow \cA$ where $\Dom(b)$ is a (possibly empty) finite subset of $\mathbb{N}_0.$ Similarly, a \textit{located  variable word } in the alphabet $\cA$ and variable $v$ is a finitely supported function $b : \Dom(b) \rightarrow \cA \cup \{v\}$ whose range contains $v$, where $\Dom(b)$ is a finite subset of $\mathbb{N}_0$. Let $L(\cA)$ be the set of \textit{located words} in $\cA$ and let $L(\cA v)$ be the set of located variable words in $\cA$ for the variable $v.$ Then $S:= L(\cA) \cup L(\cA v)$ has a natural partial semigroup operation, obtained by letting $b_0 \ast b_1$ be defined whenever the domains of $b_0$ and $b_1$ are disjoint. In such a case, $b_0 \ast b_1$ is just $b_0 \cup b_1.$

\begin{theorem}[Hales-Jewett Theorem] 
Let \( L(\mathcal{A}) \) be finitely colored. Then there exist \( \alpha \in L(\mathcal{A}) \) and \( \gamma \in \mathcal{P}_f(\mathbb{N}) \) such that \( \Dom(\alpha) \cap \gamma = \emptyset \) and \( \{\alpha \cup (\gamma \times \{s\}) : s \in \mathcal{A}\} \) is monochromatic.
\end{theorem}

In 2008, M. Beiglb\"{o}ck extended the Hales–Jewett theorem, obtaining a result stronger than the one above. His extension involves partition regular families of \( \mathbb{N} \). A \emph{partition regular family} \( \mathcal{F} \) of \( \mathbb{N} \) is a subset of \( \mathcal{P}_f(\mathbb{N}) \) such that for any partition \( \mathbb{N} = C_1 \cup \cdots \cup C_r \), there exists an \( i \in \{1, \dots, r\} \) such that \( C_i \in \mathcal{F} \).

\begin{theorem}[Beiglb\"{o}ck, 2008] \label{Beiglbock}
Let \( \mathcal{F} \) be a partition regular family of finite subsets of \( \mathbb{N} \) that contains no singletons, and let \( \mathcal{A} \) be a finite alphabet set. For any finite coloring of \( L(\mathcal{A}) \), there exist \( \alpha \in L(\mathcal{A}) \), \( \gamma \in \mathcal{P}_f(\mathbb{N}) \), and \( F \in \mathcal{F} \) such that \( \Dom(\alpha) \), \( \gamma \), and \( F \) are pairwise disjoint sets, and 
\[
\{\alpha \cup (\gamma \cup \{t\} \times \{s\}) : s \in \mathcal{A}, t \in F\}
\]
is monochromatic.
\end{theorem}

\subsection{Geo-Arithmetic Progressions}
Vitaly Bergelson proved that for any finite coloring of \( \mathbb{Z} \), there exists a monochromatic geo-arithmetic progression of arbitrary length. This result can be viewed as a combined extension of the additive and multiplicative versions of van der Waerden's theorem. Bergelson initially proved this property using ergodic theory \cite{Ber05}. Later, M. Beiglb\"{o}ck, V. Bergelson, N. Hindman, and D. Strauss provided an alternative proof using the algebra of the Stone–Čech compactification in \cite{BBHS08}.

\begin{theorem}[Geo-Arithmetic Progression] \label{Originalgeoarith}
If \( n, r \in \mathbb{N} \) and \( \mathbb{Z} \) is \( r \)-colored, then there exist \( a, b, d \in \mathbb{N} \) such that the set \( \{a \cdot (b + i \cdot d)^j : 0 \leq i, j \leq n\} \) is monochromatic.
\end{theorem}

In \cite{Beigl08}, M. Beiglb\"{o}ck demonstrated that the extension of the Hales–Jewett theorem is sufficiently strong to imply Theorem \ref{Originalgeoarith}. In this article, we will use this variant of the Hales–Jewett theorem (Theorem \ref{Beiglbock}) to prove our geo-arithmetic structure.

\begin{theorem} \label{Geoarith}
If \( k, r \in \mathbb{N} \) and \( \mathbb{N} \) is \( r \)-colored, then there exist \( a, b, d \in \mathbb{N} \) and \( \gamma \in \mathcal{P}_f(\mathbb{N}) \) such that the set 
\[
\left\{\left|\Sigma_{<s_b \left(\prod_{t \in \gamma} s_t \cdot s_{a + i d}\right)^j}\right| : i, j = 0, 1, 2, \dots, k\right\}
\]
is monochromatic, where \( \Sigma \) is the set of sums of two squares.
\end{theorem}

\begin{proof}
Assume that \( \mathbb{N} \) is finitely colored. Fix \( k \in \mathbb{N} \), and let \( \mathcal{F} = \{\{a, a + d, \dots, a + k d\} : a, d \in \mathbb{N}\} \) be the set of all \( (k + 1) \)-term arithmetic progressions. Let \( \mathcal{A} = \{0, 1, \dots, k\} \), and define the function 
\[
h : L(\mathcal{A}) \rightarrow \mathbb{N} \quad \text{by} \quad h(\alpha) = \underset{t \in \Dom(\alpha)}{\ast_f} t^{(\alpha(t))}.
\]
We color each \( \alpha \in L(\mathcal{A}) \) with the color of \( h(\alpha) \) and choose \( \alpha \), \( \gamma \), and \( F = \{a, a + d, \dots, a + k d\} \).

By Theorem \ref{Beiglbock}, the set 
\[
\{h(\alpha \cup (\gamma \cup \{a + i d\}) \times \{j\}) : i, j \in \{0, 1, \dots, k\}\}
\]
is monochromatic. Expanding \( h \), we obtain:
\[
h(\alpha \cup (\gamma \cup \{a + i d\}) \times \{j\}) = \left(\underset{t \in \Dom(\alpha)}{\ast_f} t^{(\alpha(t))}\right) \ast_f \left(\underset{t \in \gamma}{\ast_f} t^{(j)}\right) \ast_f (a + i d)^{(j)}.
\]
This can be rewritten as:
\[
\left|\Sigma_{<\left(\prod_{t \in \Dom(\alpha)} s_t^{\alpha(t)}\right) \left(\prod_{t \in \gamma} s_t \cdot s_{a + i d}\right)^j}\right|.
\]
Let \( s_b = \prod_{t \in \Dom(\alpha)} s_t^{\alpha(t)} \). Then the set 
\[
\left\{\left|\Sigma_{<s_b \left(\prod_{t \in \gamma} s_t \cdot s_{a + i d}\right)^j}\right| : i, j \in \{0, 1, \dots, k\}\right\}
\]
is monochromatic, as required.
\end{proof}

%Here one notes $t^{(\alpha(t))}= \overbrace{t \ast_f t \ast_f \cdots \ast_f t }^{\alpha(t)  \mbox{-times} }= \card(\Sigma/{s_t^{\alpha(t)}})$

%Therefore, ${\ast_f}_{t \in \Dom(\alpha)}t^{(\alpha(t))}= \card(S(2,2)/_{\prod_{t \in \Dom(\alpha)}{s_t^{\alpha(t)}}})$.

%Similarly, we have ${\ast_f}_{t \in \gamma}t^{(j)}= \card (S(2,2)/_{\prod_{t \in \gamma}{s_t^j}})$ and $(a+id)^{(j)}= \card (S(2,2)/_{s_{a+id}^j})$

%Then from \eqref{eq1}, we have $\card (S(2,2)/_{({\prod_{t \in \Dom(\alpha)}{s_t^{\alpha(t)}}} )({\prod_{t \in \gamma}{s_t^j}})( {s_{a+id}^j} )})= \card (S(2,2)/_{({\prod_{t \in \Dom(\alpha)}{s_t^{\alpha(t)}}} )({\prod_{t \in \gamma}{s_t}} {s_{a+id}} )^j})$

%The set $\{\card (S(2,2)/_{({\prod_{t \in \Dom(\alpha)}{s_t^{\alpha(t)}}} )({\prod_{t \in \gamma}{s_t}} {s_{a+id}} )^j}): i,j=0,1,2, \cdots, k\}$ is monochromatic, i.e., the set  $\{\card (S(2,2)/_{s_b({\prod_{t \in \gamma}{s_t}} {s_{a+id}} )^j}): i,j=0,1,2, \cdots, k\}$ is monochromatic, where $s_b={\prod_{t \in \Dom(\alpha)}{s_t^{\alpha(t)}}}.$

\subsection{Polynomial van der Waerden's Theorem}
The polynomial extension of van der Waerden's theorem relies on the polynomial version of the Hales-Jewett theorem. In 1988, V. Bergelson and A. Liebman proved the polynomial extension of the Hales-Jewett theorem by introducing and developing the apparatus of set-polynomials (polynomials whose coefficients are finite sets) and applying the methods of topological dynamics \cite{BerLeib89}. Later, M. Walters provided short and purely combinatorial proofs of those results in \cite{Walters2000}. Let us begin with the statement of the polynomial Hales-Jewett theorem, along with some relevant notations.

For fixed numbers \( q, N, d \), consider the set \( X(q, N, d) = \prod_{i=1}^d [q]^{N^i} \), where \( [q] = \{1, \dots, q\} \). An element \( x \in X(q, N, d) \) is of the form \( (\mathbf{b}_1, \dots, \mathbf{b}_d) \), where \( \mathbf{b}_i : [N]^i \to [q] \). For \( a = (\mathbf{a}_1, \dots, \mathbf{a}_d) \), \( \gamma \subseteq [N] \), and \( (x_1, \dots, x_d) \in [q]^d \), define an element \( x = a \oplus x_1 \gamma \oplus \cdots \oplus x_d \gamma^d \) as follows:

If \( x = (\mathbf{b}_1, \dots, \mathbf{b}_d) \), then
\[
\mathbf{b}_j((i_1, \dots, i_j)) = 
\begin{cases} 
x_j, & \text{if } (i_1, \dots, i_j) \in \gamma^j, \\
\mathbf{a}_j((i_1, \dots, i_j)), & \text{otherwise.}
\end{cases}
\]

\begin{theorem}[Polynomial Hales-Jewett (PHJ) Theorem] \label{phj}
For any \( q, k, d \in \mathbb{N} \), there exists \( N \in \mathbb{N} \) such that whenever \( X(q, N, d) \) is \( k \)-colored, there exist \( a \in X(q, N, d) \) and \( \gamma \subseteq [N] \) such that the set
\[
\{ a \oplus x_1 \gamma \oplus \cdots \oplus x_d \gamma^d : (x_1, \dots, x_d) \in [q]^d \}
\]
is monochromatic.
\end{theorem}

We will use Theorem \ref{phj} to derive our version of the polynomial van der Waerden theorem, which is stated below:

\begin{theorem}
Let \( d, k, \ell \in \mathbb{N} \) and \( \{F_1, \dots, F_\ell\} \subset \mathcal{P}_f(\mathbb{N}) \) with \( F_i = \{a_{i1}, \dots, a_{id}\} \) for all \( i \in \{1, \dots, \ell\} \). Then for any \( k \)-coloring of \( \mathbb{N} \), there exist \( b, c \in \mathbb{N} \) such that the set 
\[
\left\{\left|\Sigma_{<s_b s_{a_{i1}}^{c} \cdots s_{a_{id}}^{c^d}}\right| : i = 1, \dots, \ell \right\}
\]
is monochromatic, where \( \Sigma \) is the set of sums of two squares.
\end{theorem}

\begin{proof}
Let \( q, k, d \in \mathbb{N} \) be as in the statement of the polynomial Hales-Jewett theorem. By the theorem, there exists a natural number \( N = N(q, k, d) \). Suppose \( \chi : (\mathbb{N}, \ast_f) \to [k] \) is a \( k \)-coloring of \( \mathbb{N} \), and let \( m : X(q, N, d) \to (\mathbb{N}, \ast_f) \) be the canonical map defined by
\[
m((\mathbf{b}_1, \dots, \mathbf{b}_d)) = \overset{d}{\underset{j=1}{\ast_f}} \left( \underset{(i_1, \dots, i_j) \in [N]^j}{\ast_f} \mathbf{b}_j((i_1, \dots, i_j)) \right).
\]
The composite \( \chi \circ m \) is a \( k \)-coloring of \( X(q, N, d) \). By the polynomial Hales-Jewett theorem, there exist \( a \in X(q, N, d) \) and \( \gamma \subseteq [N] \) such that the set
\[
\{ a \oplus x_1 \gamma \oplus \cdots \oplus x_d \gamma^d : (x_1, \dots, x_d) \in [q]^d \}
\]
is monochromatic. Therefore, the image
\[
m\left(\{ a \oplus x_1 \gamma \oplus \cdots \oplus x_d \gamma^d : (x_1, \dots, x_d) \in [q]^d \}\right)
\]
is monochromatic for the coloring \( \chi \) of \( \mathbb{N} \). Note that
\[
m(a \oplus x_1 \gamma \oplus \cdots \oplus x_d \gamma^d) = \overset{d}{\underset{j=1}{\ast_f}} \left( \underset{(i_1, \dots, i_j) \in \gamma^j}{\ast_f} \mathbf{b}_j((i_1, \dots, i_j)) \right) \ast_f \left( \overset{d}{\underset{j=1}{\ast_f}} \left( \underset{(i_1, \dots, i_j) \in [N]^j \setminus \gamma^j}{\ast_f} \mathbf{b}_j((i_1, \dots, i_j)) \right) \right).
\]
This can be rewritten as
\[
x_1^{(c)} \ast_f \cdots \ast_f x_d^{(c^d)} \ast_f b = |\Sigma_{<  s_b s_{x_1}^{c} \cdots s_{x_d}^{c^d}}|,
\]
where \( c = |\gamma| \) and
\[
b = \overset{d}{\underset{j=1}{\ast_f}} \left( \underset{(i_1, \dots, i_j) \in [N]^j \setminus \gamma^j}{\ast_f} \mathbf{a}_j((i_1, \dots, i_j)) \right).
\]
\end{proof}

%\begin{remark}It is evident that \( s_n \cdot s_m = s_N \) for some \( N \). Below we provide a Python code for computing \( s_N \) and the position \( N \), which amounts to computing the inverse of \( f \) (see Definition \ref{oper}).
%\end{remark}

\begin{remark}
Let \( (P, <) \) be a countably infinite poset. Then there is a canonical bijection
\[
\phi : P \to \mathbb{N} \quad \text{given by} \quad \phi(x) = |\{ y \in P : y < x \}|,
\]
with an inverse \( \psi \) given by \( \psi(n) = (n+1)^{\text{st}} \) smallest element of \( P \). Since \( \psi \) does not preserve the multiplication of \( \mathbb{N} \), if \( P \) is a countable multiplicative subgroup of \( \mathbb{N} \), then the pullback \( \psi \) is not a semigroup homomorphism with respect to the usual multiplication. Thus, we obtain a new operation on \( \mathbb{N} \) via
\[
m \ast_\psi n = \phi(\psi(m) \cdot \psi(n))
\]
(as we did before for the case of \( \Sigma \)). Note that \( (\mathcal{P}_f(\mathbb{N}), <) \) is a countable ordered poset, and for the map \( \sigma : \mathcal{P}_f(\mathbb{N}) \to \mathbb{N} \) given by \( \sigma(F) = \sum_{n \in F} 2^n \), this map is not of the above form, although it induces operations on \( \mathcal{P}_f(\mathbb{N}) \).
\end{remark}

\begin{remark}
The set \( \Sigma = \{ a^2 + b^2 : a, b \in \mathbb{N}_0 \} \) is symmetric with respect to \( a, b \), so there is a canonical bijection between the set \( T := \{(a, b) \in \mathbb{N}_0^2 : a \leq b \} \) and \( \Sigma \). Then the map \( \psi : T \to \mathbb{N} \) given by
\[
\psi(m, 2n) = (n+1)^2 - m \quad \text{and} \quad \psi(m, 2n+1) = (n+1)(n+2) - m
\]
is a bijection that induces a new operation on \( T \). Hence, one may write many different monochromatic configurations in \( \mathbb{N} \).
\end{remark}

\section*{Acknowledgement}The first author would like to thank the Department of Mathematics, University of Haifa, for her position.

\bibliographystyle{amsplain}

\begin{thebibliography}{10}
\bibitem{BLM21} J.M. Barrett, M. Lupini, Joel Moreira, \textit{On Rado conditions for nonlinear Diophantine equations}, Eur. J. Comb. 94 (2021) 103277.
\vspace{.3cm}
\bibitem{Beigl08} M. Beiglb\"{o}ck, \textit{A variant of the Hales-Jewett theorem}, Bull. Lond. Math. Soc. 40(2008) 210-216.
\vspace{.3cm}
\bibitem{BBHS08} M. Beiglb\"{o}ck,  V. Bergelson,  N. Hindman and D. Strauss, \textit{Some new results in multiplicative and additive Ramsey theory}, Trans. Am. Math. Soc. 360(2008) 819-847.
\vspace{.3cm}
\bibitem{BerLeib89} V. Bergelson, A. Leibman, \textit{Set-polynomials and polynomial extension of the Hales-Jewett Theorem}, Ann. Math. 150(1999) 33-75.

\vspace{.3cm}
\bibitem{Ber05} V. Bergelson, \textit{Multiplicatively large sets and Ramsey theory}, Isr. J. Math. 148(2005) 23-40.

\vspace{.3cm}

\bibitem{Ber10} V. Bergelson, \textit{IP sets, dynamics, and combinatorial number theory}, in: V. Bergelson, A. Blass, M. Di Nasso, R. Jin (Eds.), Ultrafilters Across Mathematics, in: Contemp. Math., vol. 530, AMS, 2010, pp.23-47

\vspace{.3cm}
\bibitem{BHW14} V. Bergelson, N. Hindman, K. Williams, \textit{Polynomial extensions of the Milliken--Taylor theorem}, Trans. Am.Math.Soc. 366(2014) 5727-5748.

\vspace{.3cm}
\bibitem{BJM17} V. Bergelson, J.H. Johnson, J. Moreira, \textit{New polynomial and multidimensional extensions of classical partition results}, J. Comb. Theory, Ser. A 147(2017) 119-154.

\vspace{.3cm}
\bibitem{Brauer28} A. Brauer, \textit{\"{U}ber sequenzen von potenzresten}, Sitz.ber. Preuss. Akad. Wiss. Berl. Phys.-Math. KI.(1928) 9-16.

\vspace{.3cm}
\bibitem{CG24} A. Chakraborty, S. Goswami, \textit{Polynomial extension of some symmetric partition
regular structures}, Bull. Sci. math. 192(2024) 103415.

\vspace{.3cm}

\bibitem{G25} S. Goswami, \textit{Monochromatic Translated Product and Answering Sahasrabudhe’s Conjecture}, arXiv:2412.17868v1 .

\vspace{.3cm}
\bibitem{Deuber74} W. Deuber, \textit{Partitionen und lineare Gleichungssysteme}, Math.Z. 133 (1973) 109-123.

\vspace{.3cm}
\bibitem{Hind74} N. Hindman,
\textit{Finite sums from sequences within cells of a partition of $\mathbb{N}$}, J. Comb. Theory, Ser. A 17 (1974), 1-11.

\vspace{.3cm}

\bibitem{Hind11} N.Hindman, \textit{Monochromatic sums equal to products in $\bN$}, Integers 11A (2011) 10.

\vspace{.3cm}

\bibitem{HS12} N. Hindman and D. Strauss, \textit{ Algebra in the Stone-$\breve{\mbox{C}}$ech Compactification}, Theory and
    Application, second edition, de Gruyter, Berlin, 2011.

\vspace{.3cm}    

\bibitem{Lupini19} M. Lupini, \textit{Actions on semigroups and an infinitary Gowers-Hales-Jewett Ramsey theorem}, Trans. Am. Math. Soc. 371 (2019) 3083-3116.

\vspace{.3cm}
\bibitem{Joel17} J. Moreira, \textit{Monochromatic sums and products in $\bN$}, Ann. Math.185 (2017) 1069-1090.      

\vspace{.3cm}
\bibitem{Milli75} K. Milliken, \textit{Ramsey's theorem with sums or unions}, J. Comb. Theory, Ser. A, 18(1975) 276-290.

\vspace{.3cm}
\bibitem{Nasso22} Mauro Di Nasso, \textit{Infinite monochromatic patterns in the integers},  J. Comb. Theory, Ser. A 189(2022) 105610. 

\vspace{.3cm}
\bibitem{Saha18} J. Sahasrabudhe,
\textit{Exponential patterns in arithmetic Ramsey theory},
Acta Arith. 182 (2018), no. 1, 13–42.

\vspace{.3cm}

\bibitem{Tay76} A.Taylor, \textit{A canonical partition relation for finite subsets of $\omega$}, J. Comb. Theory, Ser. A 21(1976) 137-146.


\vspace{.3cm}    
\bibitem{Van27} B. L. Van der Waerden, \textit{Beweis einer baudetschen vermutung}, Nieuw Arch. Wiskd. 15 (1927), 212- 216.

\vspace{.3cm}
\bibitem{Tod10} S. Todorcevic, \textit{Introduction to Ramsey spaces}, Annals of Mathematics Studies, vol.174, Princeton University Press, Princeton, NJ, 2010.

\vspace{.3cm}
\bibitem{Walters2000} M. Walters, \textit{Combinatorial proofs of the polynomial van der Waerden Theorem and the polynomial Hales-Jewett Theorem}, J. Lond. Math. Soc. 61 (2000) 1-12.
\vspace{.5cm}
\end{thebibliography}

\end{document}